\documentclass{article}
\usepackage{amsfonts}
\usepackage{amsmath}

\newtheorem{theorem}{Theorem}

\newtheorem{note}[theorem]{Note}

\newtheorem{corollary}[theorem]{Corollary}

\newtheorem{definition}[theorem]{Definition}

\newtheorem{lemma}[theorem]{Lemma}

\newtheorem{proposition}[theorem]{Proposition}
\newtheorem{remark}[theorem]{Remark}

\newenvironment{proof}[1][Proof]{\textbf{#1.} }{\ \rule{0.5em}{0.5em}}
\numberwithin{equation}{section}

\input{tcilatex}

\begin{document}

\title{Discrete fractional Calculus and Inequalities}
\author{George A. Anastassiou \\
Department of Mathematical Sciences\\
University of Memphis\\
Memphis, TN 38152, U.S.A.\\
ganastss@memphis.edu}
\date{}
\maketitle

\begin{abstract}
Here we define a Caputo like discrete fractional difference and we compare
it to the earlier defined Riemann-Liouville fractional discrete analog. Then
we produce discrete fractional Taylor formulae for the first time, and we
estimate their remainders. Finally, we derive related discrete fractional
Ostrowski, Poincare and Sobolev type inequalities.
\end{abstract}



\section{Preliminaries}

We make

\begin{definition}
We follow \cite{2}, \cite{3}, \cite{4}.

Let $\nu >0$. The $\nu $-th fractional sum of $f$ is defined by 
\begin{equation*}
\Delta ^{-\nu }f\left( t,a\right) =\frac{1}{\Gamma \left( \nu \right) }%
\sum_{s=a}^{t-\nu }\left( t-s-1\right) ^{\left( \nu -1\right) }f\left(
s\right) .
\end{equation*}%
Here $f$ is defined for $s=a$ $\func{mod}$ $(1)$ and $\Delta ^{-\nu }f$ is
defined for $t=\left( a+\nu \right) $ $\func{mod}$ $\left( 1\right) $; in
particular $\Delta ^{-\nu }$ maps functions defined on $\mathbb{N}_{\alpha }$
to functions defined on $\mathbb{N}_{\alpha +\nu }$, where $\mathbb{N}%
_{t}=\{t,t+1,t+2,...\}.$

Here $t^{\left( \nu \right) }=\frac{\Gamma \left( t+1\right) }{\Gamma \left(
t-\nu +1\right) }.$
\end{definition}

From now in this context for convinience we set $\Delta ^{-\nu }f\left(
t,a\right) =\Delta ^{-\nu }f\left( t\right) $.

We mention

\begin{theorem}
\label{t2}(\cite{2}) Let $f$ be a real-valued function defined on $\mathbb{N}%
_{a}$ and let $\mu ,\nu >0$. Then 
\begin{equation*}
\Delta ^{-\nu }\left( \Delta ^{-\mu }f\left( t\right) \right) =\Delta
^{-\left( \mu +\nu \right) }f\left( t\right) =\Delta ^{-\mu }\left( \Delta
^{-\nu }f\left( t\right) \right) ,\text{ \ \ \ \ \ \ \ \ \ }\forall \text{ }%
t\in \mathbb{N}_{a+\mu +\nu }.
\end{equation*}
\end{theorem}

We make

\begin{definition}
\label{d3}Let $\mu >0$ and $m-1<\mu <m$, where $m$ denotes a positive
integer, $m=\left\lceil \mu \right\rceil $, $\left\lceil .\right\rceil $
ceiling of number. Set $\nu =m-\mu $.

The $\mu $-th fractional Caputo like difference is defined as 
\begin{equation*}
\Delta _{\ast }^{\mu }f\left( t\right) =\Delta ^{-\nu }\left( \Delta
^{m}f\left( t\right) \right) =\frac{1}{\Gamma \left( \nu \right) }%
\sum_{s=a}^{t-\nu }\left( t-s-1\right) ^{\left( \nu -1\right) }\left( \Delta
^{m}f\right) \left( s\right) \text{, \ \ }\forall \text{ }t\in \mathbb{N}%
_{a+\nu }.
\end{equation*}%
Here $\Delta ^{m}$ is the $m$-th order formward difference operator 
\begin{equation*}
\left( \Delta ^{m}f\right) \left( s\right) =\sum_{k=0}^{m}\left( 
\begin{array}{c}
m \\ 
k%
\end{array}%
\right) \left( -1\right) ^{m-k}f\left( s+k\right) .
\end{equation*}
\end{definition}

We mention

\begin{theorem}
\label{t4}(\cite{3}) For $\nu >0$ and $p$ a positive integer we have 
\begin{equation*}
\Delta ^{-\nu }\Delta ^{p}f\left( t\right) =\Delta ^{p}\Delta ^{-\nu
}f\left( t\right) -\sum_{k=0}^{\nu -1}\frac{\left( t-a\right) ^{\left( \nu
-p+k\right) }}{\Gamma \left( \nu +k-p+1\right) }\Delta ^{k}f\left( a\right) 
\text{,}
\end{equation*}%
where $f$ is defined on $\mathbb{N}_{a}$.
\end{theorem}

\begin{remark}
Let $\mu >0$ and $m-1<\mu <m$, $m=\left\lceil \mu \right\rceil $, where $m$
is a positive integer, $\nu =m-\mu >0$. Then by Theorem \ref{t4} we get 
\begin{equation*}
\Delta ^{-\nu }\Delta ^{m}f\left( t\right) =\Delta ^{m}\Delta ^{-\nu
}f\left( t\right) -\sum_{k=0}^{m-1}\frac{\left( t-a\right) ^{\left( \nu
-m+k\right) }}{\Gamma \left( \nu +k-m+1\right) }\Delta ^{k}f\left( a\right) ,
\end{equation*}%
where $f$ is defined on $\mathbb{N}_{a}$.

So we have proved 
\begin{equation*}
\Delta _{\ast }^{\mu }f\left( t\right) =\Delta ^{m}\Delta ^{-\nu }f\left(
t\right) -\sum_{k=0}^{m-1}\frac{\left( t-a\right) ^{\left( \nu -m+k\right) }%
}{\Gamma \left( \nu +k-m+1\right) }\Delta ^{k}f\left( a\right) ,
\end{equation*}%
that is 
\begin{equation}
\Delta ^{m}\Delta ^{-\nu }f\left( t\right) =\Delta _{\ast }^{\mu }f\left(
t\right) +\sum_{k=0}^{m-1}\frac{\left( t-a\right) ^{\left( \nu -m+k\right) }%
}{\Gamma \left( \nu +k-m+1\right) }\Delta ^{k}f\left( a\right) ,  \tag{1}
\label{r1}
\end{equation}%
where $f$ is defined on $\mathbb{N}_{a}$.
\end{remark}

\begin{definition}
(\cite{3}) The $\mu $-th fractional Riemann-Liouville type difference is
defined by 
\begin{equation*}
\Delta ^{\mu }f\left( t\right) :=\Delta ^{m-\nu }f\left( t\right) :=\Delta
^{m}\left( \Delta ^{-\nu }f\left( t\right) \right) ,
\end{equation*}%
where $\mu >0$, $m-1<\mu <m$, $\nu =m-\mu >0$.
\end{definition}

\begin{remark}
Consequently from (\ref{r1}) we get 
\begin{equation}
\Delta ^{\mu }f\left( t\right) =\Delta _{\ast }^{\mu }f\left( t\right)
+\sum_{k=0}^{m-1}\frac{\left( t-a\right) ^{\left( \nu -m+k\right) }}{\Gamma
\left( \nu +k-m+1\right) }\Delta ^{k}f\left( a\right) ,  \tag{2}  \label{r2}
\end{equation}%
where $f$ is defined on $\mathbb{N}_{a}$.
\end{remark}

\section{Results}

We give the following Caputo type fractional Taylor's difference formula.

\begin{theorem}
\label{t8}For $\mu >0$, $\mu $ non-integer, $m=\left\lceil \mu \right\rceil $%
, $\nu =m-\mu $, it holds: 
\begin{equation}
f\left( t\right) =\sum_{k=0}^{m-1}\frac{\left( t-a\right) ^{\left( k\right) }%
}{k!}\Delta ^{k}f\left( a\right) +\frac{1}{\Gamma \left( \mu \right) }%
\sum_{s=a+\nu }^{t-\mu }\left( t-s-1\right) ^{\left( \mu -1\right) }\Delta
_{\ast }^{\mu }f\left( s\right) ,\text{ }\forall t\in \mathbb{N}_{a+m}, 
\tag{3}  \label{r3}
\end{equation}%
where $f$ is defined on $\mathbb{N}_{a}$ with $a\in \mathbb{Z}^{+}$, $%
\mathbb{Z}^{+}:=\{0,1,2,...\}.$
\end{theorem}

\begin{proof}
Notice that by Definition \ref{d3},

\noindent $\Delta _{\ast }^{\mu }f\left( t\right) =\Delta ^{-\nu }\left(
\Delta ^{m}f\left( t\right) \right) =\Delta ^{-\left( m-\mu \right) }\left(
\Delta ^{m}f\left( t\right) \right) $, $\forall $ $t\in \mathbb{N}_{a+\nu }$.

Consequently we get $\Delta ^{-\mu }\Delta _{\ast }^{\mu }f\left( t\right)
=\Delta ^{-\mu }\Delta ^{-\left( m-\mu \right) }\left( \Delta ^{m}f\left(
t\right) \right) $ (by Theorem \ref{t2}) $=\Delta ^{-\left( \mu +\left(
m-\mu \right) \right) }\left( \Delta ^{m}f\left( t\right) \right) =\Delta
^{-m}\left( \Delta ^{m}f\left( t\right) \right) $, $\forall $ $t\in \mathbb{N%
}_{a+\nu +\mu }$.

So that 
\begin{equation}
\Delta ^{-\mu }\Delta _{\ast }^{\mu }f\left( t\right) =\Delta ^{-m}\left(
\Delta ^{m}f\left( t\right) \right) ,\text{ \ \ \ \ }\forall \text{ }t\in 
\mathbb{N}_{a+m}.  \tag{4}  \label{r4}
\end{equation}%
We notice that 
\begin{equation}
\left( t-s-1\right) ^{\left( m-1\right) }=\frac{\Gamma \left( t-s\right) }{%
\Gamma \left( t-s-m+1\right) }=\left( t-s-1\right) \left( t-s-2\right)
...\left( t-s-m+1\right) ,  \tag{5}  \label{r5}
\end{equation}%
the falling factorial, here we have $t-s-m+1>0.$

Therefore we obtain 
\begin{equation}
\Delta ^{-m}\left( \Delta ^{m}f\left( t\right) \right) =\frac{1}{\left(
m-1\right) !}\sum_{s=a}^{t-m}\left( t-s-1\right) ^{\left( m-1\right) }\Delta
^{m}f\left( s\right) .  \tag{6}  \label{r6}
\end{equation}%
By (\cite{1}, p. 28, Theorem 1.8.5) the discrete Taylor's formula we get 
\begin{equation}
f\left( t\right) =\sum_{k=0}^{m-1}\frac{\left( t-a\right) ^{\left( k\right) }%
}{k!}\Delta ^{k}f\left( a\right) +\frac{1}{\left( m-1\right) !}%
\sum_{s=a}^{t-m}\left( t-s-1\right) ^{\left( m-1\right) }\Delta ^{m}f\left(
s\right) ,  \tag{7}  \label{r7}
\end{equation}%
where $t^{\left( k\right) }=t\left( t-1\right) ...\left( t-k+1\right) .$

From the last we derive 
\begin{equation}
f\left( t\right) =\sum_{k=0}^{m-1}\frac{\left( t-a\right) ^{\left( k\right) }%
}{k!}\Delta ^{k}f\left( a\right) +\Delta ^{-\mu }\Delta _{\ast }^{\mu
}f\left( t\right) ,  \tag{8}  \label{r8}
\end{equation}%
where $f$ is defined on $\mathbb{N}_{a}$, $\forall $ $t\in \mathbb{N}_{a+m},$
proving the claim.
\end{proof}

We make

\begin{remark}
Here $[a,b]$ denotes the discrete interval $[a,b]=[a,a+1,a+2,...,b]$, where $%
a<b$ and $a,b\in \{0,1,...\}$.

Let $\mu >0$ be non integer such that $m-1<\mu <m$, i.e. $m=\left\lceil \mu
\right\rceil $. Consider a function $f$ defined on $[a,b]$. Then clearly the
fractional discrete Taylor's formula (\ref{r3}) is valid only for $t\in
\lbrack a+m,b]$, $a+m<b$.

We need
\end{remark}

\begin{theorem}
(\cite{3}) Let $p$ be a positive integer and let $\nu >p$. Then 
\begin{equation}
\Delta ^{p}\left( \Delta ^{-\nu }f\left( t\right) \right) =\Delta ^{-\left(
\nu -p\right) }f\left( t\right) .  \tag{9}  \label{r9}
\end{equation}
\end{theorem}

We make

\begin{remark}
Let $\mu >p$, where $p\in \mathbb{N}$. Then 
\begin{equation}
\Delta ^{p}\left( \Delta ^{-\mu }\Delta _{\ast }^{\mu }f\left( t\right)
\right) \overset{(\ref{r9})}{=}\Delta ^{-\left( \mu -p\right) }\left( \Delta
_{\ast }^{\mu }f\left( t\right) \right) ,\text{ \ \ \ \ }\forall \text{ }%
t\in \mathbb{N}_{a+m-p}.  \tag{10}  \label{r10}
\end{equation}%
Also notice that 
\begin{equation}
\Delta ^{p}\left( \frac{\left( t-a\right) ^{\left( k\right) }}{k!}\right) =%
\frac{\left( t-a\right) ^{\left( k-p\right) }}{\left( k-p\right) !}\text{, \
\ \ \ for }k\geq p.  \tag{11}  \label{r11}
\end{equation}%
By the last we obtain the following discrete Caputo type fractional extended
Taylor's formula.
\end{remark}

\begin{theorem}
\label{t12}Let $\mu >p$, $p\in \mathbb{N}$, $\mu $ not integer, $%
m=\left\lceil \mu \right\rceil $, $\nu =m-\mu $. Then 
\begin{equation}
\Delta ^{p}f\left( t\right) =\sum_{k=p}^{m-1}\frac{\left( t-a\right)
^{\left( k-p\right) }}{\left( k-p\right) !}\Delta ^{k}f\left( a\right) +%
\frac{1}{\Gamma \left( \mu -p\right) }\sum_{s=a+\nu }^{t-\mu +p}\left(
t-s-1\right) ^{\left( \mu -p-1\right) }\Delta _{\ast }^{\mu }f\left(
s\right) ,  \tag{12}  \label{r12}
\end{equation}%
$\forall $ $t\in \mathbb{N}_{a+m-p}$, $f$ is defined on $\mathbb{N}_{a}$, $%
a\in \mathbb{Z}^{+}$.
\end{theorem}

\begin{note}
Assuming that $f$ is defined on $[a,b]$, then (\ref{r12}) is valid only for $%
[a+m-p,b]$, with $a+m-p<b$.

Notice for $p=0$ applied on (\ref{r12}) we get (\ref{r3}).
\end{note}

We give

\begin{proposition}
For $\mu >0$, $\mu $ not an integer, $m=\left\lceil \mu \right\rceil $, $\nu
=m-\mu $, $f$ is defined on $\mathbb{N}_{a}$, $a\in \mathbb{Z}^{+}$; and $%
\Delta ^{k}f\left( a\right) =0$, for $k=0,...,m-1,$ we get 
\begin{equation}
f\left( t\right) =\frac{1}{\Gamma \left( \mu \right) }\sum_{s=a+\nu }^{t-\mu
}\left( t-s-1\right) ^{\left( \mu -1\right) }\Delta _{\ast }^{\mu }f\left(
s\right) ,\text{ \ \ \ \ }\forall \text{ }t\in \mathbb{N}_{a+m}.  \tag{13}
\label{r13}
\end{equation}
\end{proposition}

\begin{proof}
By (\ref{r3}).
\end{proof}

Also we present

\begin{proposition}
Let $\mu >p$, $p\in \mathbb{N}$, $\mu $ non-integer, $m=\left\lceil \mu
\right\rceil $, $\nu =m-\mu $; $f$ is defined on $\mathbb{N}_{a}$, $a\in 
\mathbb{Z}^{+}$. Assume that $\Delta ^{k}f\left( a\right) =0$, $k=p,...,m-1$%
. Then 
\begin{equation}
\Delta ^{p}f\left( t\right) =\frac{1}{\Gamma \left( \mu -p\right) }%
\sum_{s=a+\nu }^{t-\mu +p}\left( t-s-1\right) ^{\left( \mu -p-1\right)
}\Delta _{\ast }^{\mu }f\left( s\right) ,\text{ \ \ \ \ }\forall \text{ }%
t\in \mathbb{N}_{a+m-p}.  \tag{14}  \label{r14}
\end{equation}
\end{proposition}

\begin{proof}
By (\ref{r12}).
\end{proof}

We make

\begin{remark}
We want to find 
\begin{equation}
\sum_{s=a+\nu }^{t-\mu }\left( t-s-1\right) ^{\left( \mu -1\right) }\text{=}%
\sum_{s=a+\nu }^{t-\mu }\frac{\Gamma \left( t-s\right) }{\Gamma \left(
t-s-\mu +1\right) }\text{=}\sum_{s=a+\nu }^{t-\mu -1}\frac{\Gamma \left(
t-s\right) }{\Gamma \left( t-s-\mu +1\right) }+\Gamma \left( \mu \right) . 
\tag{15}  \label{r15}
\end{equation}%
We notive that 
\begin{equation}
\frac{\Gamma \left( x+1\right) }{\Gamma \left( k+1\right) \Gamma \left(
x-k+1\right) }=\frac{\Gamma \left( x+2\right) }{\Gamma \left( k+2\right)
\Gamma \left( x-k+1\right) }-\frac{\Gamma \left( x+1\right) }{\Gamma \left(
k+2\right) \Gamma \left( x-k\right) }  \tag{16}  \label{r16}
\end{equation}%
with $x>k$, $x,k\in \mathbb{R}$; $k>-1$, $x>-1.$

That is 
\begin{equation}
\frac{\Gamma \left( x+1\right) }{\Gamma \left( x-k+1\right) }=\frac{1}{%
\left( k+1\right) }\left( \frac{\Gamma \left( x+2\right) }{\Gamma \left(
x-k+1\right) }-\frac{\Gamma \left( x+1\right) }{\Gamma \left( x-k\right) }%
\right) .  \tag{17}  \label{r17}
\end{equation}%
We find $A:=\sum_{s=a+\nu }^{t-\mu -1}\frac{\Gamma \left( t-s\right) }{%
\Gamma \left( t-s-\mu +1\right) }=$

\noindent (by (\ref{r17}) for $x:=t-s-1\geq \mu >0$, $k:=\mu -1>-1$, and $%
x>k $)

\noindent $\frac{1}{\mu }\sum_{s=a+\nu }^{t-\mu -1}\left[ \frac{\Gamma
\left( t-s+1\right) }{\Gamma \left( t-s+1-\mu \right) }-\frac{\Gamma \left(
t-s\right) }{\Gamma \left( t-s-\mu \right) }\right] =$

\noindent $\frac{1}{\mu }\left[ \left( \frac{\Gamma \left( t-a-\nu +1\right) 
}{\Gamma \left( t-a-\nu +1-\mu \right) }-\frac{\Gamma \left( t-a-\nu \right) 
}{\Gamma \left( t-a-\nu -\mu \right) }\right) \right. +\left( \frac{\Gamma
\left( t-a-\nu \right) }{\Gamma \left( t-a-\nu -\mu \right) }-\frac{\Gamma
\left( t-a-\nu -1\right) }{\Gamma \left( t-a-\nu -1-\mu \right) }\right) +$

\noindent $\left( \frac{\Gamma \left( t-a-\nu -1\right) }{\Gamma \left(
t-a-\nu -1-\mu \right) }-\frac{\Gamma \left( t-a-\nu -2\right) }{\Gamma
\left( t-a-\nu -2-\mu \right) }\right) +...\left. \left( ...-\frac{\Gamma
\left( \mu +1\right) }{\Gamma \left( 1\right) }\right) \right] $ $=$

\noindent $\left[ \frac{\Gamma \left( t-a-\nu +1\right) }{\mu \Gamma \left(
t-a-\nu +1-\mu \right) }-\Gamma \left( \mu \right) \right] $.

That is 
\begin{equation}
A=\frac{\Gamma \left( t-a-\nu +1\right) }{\mu \Gamma \left( t-a-\nu +1-\mu
\right) }-\Gamma \left( \mu \right) .  \tag{18}  \label{r18}
\end{equation}

Consequently we found 
\begin{equation}
\sum_{s=a+\nu }^{t-\mu }\left( t-s-1\right) ^{\left( \mu -1\right) }=\frac{%
\Gamma \left( t-a-\nu +1\right) }{\mu \Gamma \left( t-a+1-m\right) }=\frac{%
\left( t-a-\nu \right) ^{(\mu )}}{\mu }.  \tag{19}  \label{r19}
\end{equation}
\end{remark}

Using (\ref{r19}) we give

\begin{corollary}
(to Theorem \ref{t8}) Let $\mu >0$, $\mu $ non-integer, $m=\left\lceil \mu
\right\rceil $, $\nu =m-\mu $, $t\in \mathbb{N}_{a+m}$, $f$ defined on $%
\mathbb{N}_{a}$, $a\in \mathbb{Z}^{+}$. Then 
\begin{equation}
\left\vert f\left( t\right) -\sum_{k=0}^{m-1}\frac{\left( t-a\right)
^{\left( k\right) }}{k!}\Delta ^{k}f\left( a\right) \right\vert \leq \frac{%
\left( t-a-\nu \right) ^{(\mu )}}{\Gamma \left( \mu +1\right) }\cdot 
\underset{s\in \{a+\nu ,a+\nu +1,...,t-\mu \}}{\max }\left\vert \Delta
_{\ast }^{\mu }f\left( s\right) \right\vert .  \tag{20}  \label{r20}
\end{equation}
\end{corollary}

Similarly we get

\begin{corollary}
(to Theorem \ref{t12}) Let $\mu >p$, $p\in \mathbb{N}$, $\mu $ non-integer, $%
m=\left\lceil \mu \right\rceil $, $\nu =m-\mu $, $t\in \mathbb{N}_{a+m-p}$, $%
f$ defined on $\mathbb{N}_{a}$, $a\in \mathbb{Z}^{+}$. Then 
\begin{equation}
\left\vert \Delta ^{p}f\left( t\right) -\sum_{k=p}^{m-1}\frac{\left(
t-a\right) ^{\left( k-p\right) }}{\left( k-p\right) !}\Delta ^{k}f\left(
a\right) \right\vert \leq \frac{\left( t-a-\nu \right) ^{(\mu -p)}}{\Gamma
\left( \mu -p+1\right) }\cdot \underset{s\in \{a+\nu ,...,t-\mu +p\}}{\max }%
\left\vert \Delta _{\ast }^{\mu }f\left( s\right) \right\vert .  \tag{21}
\label{r21}
\end{equation}
\end{corollary}

We need

\begin{lemma}
Let $a>\nu $, $a,\nu >-1$, $a,\nu \in \mathbb{R}$, $a\leq b$. Then 
\begin{equation}
\sum_{r=a}^{b}r^{\left( \nu \right) }=\frac{1}{\left( \nu +1\right) }\left( 
\frac{\Gamma \left( b+2\right) }{\Gamma \left( b-\nu +1\right) }-\frac{%
\Gamma \left( a+1\right) }{\Gamma \left( a-\nu \right) }\right) =\left( 
\frac{\left( b+1\right) ^{\left( \nu +1\right) }-a^{\left( \nu +1\right) }}{%
\nu +1}\right) .  \tag{22}  \label{r22}
\end{equation}
\end{lemma}

\begin{proof}
We have

\noindent $\sum_{r=a}^{b}r^{\left( \nu \right) }=\sum_{r=a}^{b}\frac{\Gamma
\left( r+1\right) }{\Gamma \left( r-\nu +1\right) }\overset{\text{(by (\ref%
{r17}))}}{=}$ $\frac{1}{\left( \nu +1\right) }\sum_{r=a}^{b}\left( \frac{%
\Gamma \left( r+2\right) }{\Gamma \left( r-\nu +1\right) }-\frac{\Gamma
\left( r+1\right) }{\Gamma \left( r-\nu \right) }\right) =$

\noindent $\frac{1}{\left( \nu +1\right) }\left( \sum_{r=a}^{b}\left\{
\left( \frac{\Gamma \left( a+2\right) }{\Gamma \left( a-\nu +1\right) }-%
\frac{\Gamma \left( a+1\right) }{\Gamma \left( a-\nu \right) }\right)
\right. \right. +\left( \frac{\Gamma \left( a+3\right) }{\Gamma \left(
a+2-\nu \right) }-\frac{\Gamma \left( a+2\right) }{\Gamma \left( a+1-\nu
\right) }\right) +$

\noindent $\left( \frac{\Gamma \left( a+4\right) }{\Gamma \left( a+3-\nu
\right) }-\frac{\Gamma \left( a+3\right) }{\Gamma \left( a+2-\nu \right) }%
\right) +...+\left( \frac{\Gamma \left( b+1\right) }{\Gamma \left( b-\nu
\right) }-\frac{\Gamma \left( b\right) }{\Gamma \left( b-1-\nu \right) }%
\right) +\left. \left. \left( \frac{\Gamma \left( b+2\right) }{\Gamma (b-\nu
+1)}-\frac{\Gamma \left( b+1\right) }{\Gamma \left( b-\nu \right) }\right)
\right\} \right) =$

\noindent $\frac{1}{\left( \nu +1\right) }\left( \frac{\Gamma \left(
b+2\right) }{\Gamma \left( b-\nu +1\right) }-\frac{\Gamma \left( a+1\right) 
}{\Gamma \left( a-\nu \right) }\right) $,

\noindent proving the claim.
\end{proof}

Next we present a discrete fractional Ostrowski type inequality.

\begin{theorem}
Let $\mu >p$, $p\in \mathbb{Z}^{+}$, $\mu $ not an integer, $m=\left\lceil
\mu \right\rceil $, $\nu =m-\mu $. Here $f$ is defined on $\mathbb{N}_{a}$, $%
a\in \mathbb{Z}^{+}$ and $j\in \lbrack a+m-p+1,b]$, with $a+m-p<b\in \mathbb{%
N}$. Assume that $\Delta ^{k}f\left( a\right) =0$, for $k\in \lbrack
p+1,...,m-1]$.

Then 
\begin{equation*}
\left| \left( \frac{1}{\left( b-a-m+p\right) }\sum_{j=a+m-p+1}^{b}\Delta
^{p}f\left( j\right) \right) -\Delta ^{p}f\left( a\right) \right| \leq
\end{equation*}%
\begin{equation*}
\frac{1}{\left( b-a-m+p\right) \Gamma \left( \mu -p+2\right) }\left[ \left(
b-a-\nu +1\right) ^{\left( \mu -p+1\right) }-\Gamma \left( \mu -p+2\right) %
\right]
\end{equation*}%
\begin{equation}
\cdot \left( \underset{t\in \{a+\nu ,...,b-\mu +p\}}{\max }\left| \Delta
_{\ast }^{\mu }f\left( t\right) \right| \right) .  \tag{23}  \label{r23}
\end{equation}
\end{theorem}

\begin{proof}
By (\ref{r12}) we have 
\begin{equation}
\Delta ^{p}f\left( j\right) -\Delta ^{p}f\left( a\right) =\frac{1}{\Gamma
\left( \mu -p\right) }\sum_{s=a+\nu }^{j-\mu +p}\left( j-s-1\right) ^{\left(
\mu -p-1\right) }\Delta _{\ast }^{\mu }f\left( s\right) ,  \tag{24}
\label{r24}
\end{equation}%
for all $j\in \lbrack a+m-p+1,b].$

We get that 
\begin{eqnarray*}
\frac{1}{b-\left( a+m-p\right) }\sum_{j=a+m-p+1}^{b}\Delta ^{p}f\left(
j\right) -\Delta ^{p}f\left( a\right) &=& \\
\frac{1}{b-\left( a+m-p\right) }\sum_{j=a+m-p+1}^{b}\left( \Delta
^{p}f\left( j\right) -\Delta ^{p}f\left( a\right) \right) &=&
\end{eqnarray*}%
\begin{equation}
\frac{1}{\left( b-\left( a+m-p\right) \right) \Gamma \left( \mu -p\right) }%
\sum_{j=a+m-p+1}^{b}\left( \sum_{s=a+\nu }^{j-\mu +p}\left( j-s-1\right)
^{\left( ^{\mu -p-1}\right) }\Delta _{\ast }^{\mu }f\left( s\right) \right) .
\tag{25}  \label{r25}
\end{equation}%
Therefore we obtain

\noindent $\left| \frac{1}{b-\left( a+m-p\right) }\sum_{j=a+m-p+1}^{b}\Delta
^{p}f\left( j\right) -\Delta ^{p}f\left( a\right) \right| =$

\noindent $\frac{1}{\left( b-a-m+p\right) }\left| \sum_{j=a+m-p+1}^{b}\left(
\Delta ^{p}f\left( j\right) -\Delta ^{p}f\left( a\right) \right) \right|
\leq $

\noindent $\frac{1}{\left( b-a-m+p\right) }\sum_{j=a+m-p+1}^{b}\left| \Delta
^{p}f\left( j\right) -\Delta ^{p}f\left( a\right) \right| \leq $

\noindent $\frac{1}{\left( b-a-m+p\right) \Gamma \left( \mu -p\right) }%
\sum_{j=a+m-p+1}^{b}\left( \sum_{s=a+\nu }^{j-\mu +p}\left( j-s-1\right)
^{\left( \mu -p-1\right) }\cdot \left| \Delta _{\ast }^{\mu }f\left(
s\right) \right| \right) \overset{\text{(by (\ref{r19}))}}{\leq }$

\noindent $\frac{1}{\left( b-a-m+p\right) \Gamma \left( \mu -p+1\right) }%
\sum_{j=a+m-p+1}^{b}\left( j-a-\nu \right) ^{\left( \mu -p\right) }\cdot 
\underset{s\in \{a+\nu ,...j-\mu +p\}}{\max }\left\vert \Delta _{\ast }^{\mu
}f\left( s\right) \right\vert \leq $%
\begin{equation}
\frac{1}{\left( b\text{-}a\text{-}m\text{+}p\right) \Gamma \left( \mu \text{-%
}p\text{+}1\right) }\left( \sum_{j=a+m-p+1}^{b}\left( j\text{-}a\text{-}\nu
\right) ^{\left( \mu \text{-}p\right) }\right) \cdot \underset{s\in \{a\text{%
+}\nu ,...,b\text{-}\mu \text{+}p\}}{\max }\left\vert \Delta _{\ast }^{\mu
}f\left( s\right) \right\vert =:\left( \ast \right)  \tag{26}  \label{r26}
\end{equation}%
Next we use (\ref{r22}). We notice that 
\begin{equation*}
\sum_{j=a+m-p+1}^{b}\left( j-a-\nu \right) ^{\left( \mu -p\right)
}=\sum_{r=\mu -p+1}^{b-a-\nu }r^{\left( \mu -p\right) }=
\end{equation*}%
\begin{equation}
\frac{1}{\left( \mu -p+1\right) }\left( \frac{\Gamma \left( b-a-\nu
+2\right) }{\Gamma \left( b-a-m+p+1\right) }-\Gamma \left( \mu -p+2\right)
\right) .  \tag{27}  \label{r27}
\end{equation}%
Therefore 
\begin{equation*}
\left( \ast \right) =\frac{1}{\left( b-a-m+p\right) \Gamma \left( \mu
-p+2\right) }\left( \frac{\Gamma \left( b-a-\nu +2\right) }{\Gamma \left(
b-a-m+p+1\right) }-\Gamma \left( \mu -p+2\right) \right)
\end{equation*}%
\begin{equation}
\cdot \left( \underset{t\in \{a+\nu ,...,b-\mu +p\}}{\max }\left\vert \Delta
_{\ast }^{\mu }f\left( t\right) \right\vert \right) .  \tag{28}  \label{r28}
\end{equation}

The last completes the proof.
\end{proof}

Next we give a discrete fractional Poincar\'{e} inequality.

\begin{theorem}
Let $\mu >p$, $p\in \mathbb{Z}^{\mathbb{+}}$, $\mu $ non-integer, $%
m=\left\lceil \mu \right\rceil $, $\nu =m-\mu $. Assume that $\Delta
^{k}f\left( a\right) =0$, $k=p,...,m-1$, $f$ defined on $\mathbb{N}_{a}$, $%
a\in \mathbb{Z}^{+}$. Let $\gamma ,\delta >1:\frac{1}{\gamma }+\frac{1}{%
\delta }=1$. Then 
\begin{equation*}
\sum_{j=a+m-p}^{b}\left| \Delta ^{p}f\left( j\right) \right| ^{\delta }\leq 
\frac{1}{\left( \Gamma \left( \mu -p\right) \right) ^{\delta }}\left[
\sum_{j=a+m-p}^{b}\left( \sum_{s=a+\nu }^{j-\mu +p}\left( \left(
j-s-1\right) ^{\left( \mu -p-1\right) }\right) ^{\gamma }\right) ^{\frac{%
\delta }{\gamma }}\right]
\end{equation*}%
\begin{equation}
\cdot \left( \sum_{s=a+\nu }^{b-\mu +p}\left| \Delta _{\ast }^{\mu }f\left(
s\right) \right| ^{\delta }\right) .  \tag{29}  \label{r29}
\end{equation}
\end{theorem}

\begin{proof}
We have 
\begin{equation}
\Delta ^{p}f\left( j\right) =\frac{1}{\Gamma \left( \mu -p\right) }%
\sum_{s=a+\nu }^{j-\mu +p}\left( j-s-1\right) ^{\left( \mu -p-1\right)
}\Delta _{\ast }^{\mu }f\left( s\right) ,\text{ \ \ }\forall \text{ }j\in
\lbrack a+m-p,b].  \tag{30}  \label{r30}
\end{equation}

Let $\gamma ,\delta >1$ such that $\frac{1}{\gamma }+\frac{1}{\delta }=1.$

We observe that 
\begin{equation*}
\left| \Delta ^{p}f\left( j\right) \right| \leq \frac{1}{\Gamma \left( \mu
-p\right) }\sum_{s=a+\nu }^{j-\mu +p}\left( j-s-1\right) ^{\left( \mu
-p-1\right) }\left| \Delta _{\ast }^{\mu }f\left( s\right) \right|
\end{equation*}%
(by discrete H\"{o}lder's inequality) 
\begin{equation}
\leq \frac{1}{\Gamma \left( \mu -p\right) }\left( \sum_{s=a+\nu }^{j-\mu
+p}\left( \left( j-s-1\right) ^{\left( \mu -p-1\right) }\right) ^{\gamma
}\right) ^{\frac{1}{\gamma }}\cdot \left( \sum_{s=a+\nu }^{j-\mu +p}\left|
\Delta _{\ast }^{\mu }f\left( s\right) \right| ^{\delta }\right) ^{\frac{1}{%
\delta }}.  \tag{31}  \label{r31}
\end{equation}%
I.e. it holds 
\begin{equation*}
\left| \Delta ^{p}f\left( j\right) \right| ^{\delta }\leq \frac{1}{\left(
\Gamma \left( \mu -p\right) \right) ^{\delta }}\left( \sum_{s=a+\nu }^{j-\mu
+p}\left( \left( j-s-1\right) ^{\left( \mu -p-1\right) }\right) ^{\gamma
}\right) ^{\frac{\delta }{\gamma }}\cdot \left( \sum_{s=a+\nu }^{j-\mu
+p}\left| \Delta _{\ast }^{\mu }f\left( s\right) \right| ^{\delta }\right)
\end{equation*}%
\begin{equation*}
\leq \frac{1}{\left( \Gamma \left( \mu -p\right) \right) ^{\delta }}\left(
\sum_{s=a+\nu }^{j-\mu +p}\left( \left( j-s-1\right) ^{\left( \mu
-p-1\right) }\right) ^{\gamma }\right) ^{\frac{\delta }{\gamma }}
\end{equation*}%
\begin{equation}
\cdot \left( \sum_{s=a+\nu }^{b-\mu +p}\left| \Delta _{\ast }^{\mu }f\left(
s\right) \right| ^{\delta }\right) \text{, \ \ }\forall \text{ }j\in \lbrack
a+m-p,b].  \tag{32}  \label{r32}
\end{equation}%
Applying $\sum_{j=a+m-p}^{b}$ on (\ref{r32}) we establish (\ref{r29}).
\end{proof}

It follows a discrete Sobolev type fractional inequality.

\begin{theorem}
\label{t22}Let $\mu >p$, $p\in \mathbb{Z}^{+}$, $\mu $ non-integer, $%
m=\left\lceil \mu \right\rceil $, $\nu =m-\mu $. Assume that $\Delta
^{k}f\left( a\right) =0$, $k=p,...,m-1$; $f$ defined on $\mathbb{N}_{a}$, $%
a\in \mathbb{Z}^{+}$. Let $\gamma ,\delta >1:\frac{1}{\gamma }+\frac{1}{%
\delta }=1$, and $r\geq 1$. Then 
\begin{equation*}
\left( \sum_{j=a+m-p}^{b}\left| \Delta ^{p}f\left( j\right) \right|
^{r}\right) ^{\frac{1}{r}}\leq
\end{equation*}%
\begin{equation}
\frac{1}{\Gamma \left( \mu -p\right) }\left[ \sum_{j=a+m-p}^{b}\left(
\sum_{s=a+\nu }^{j-\mu +p}\left( \left( j-s-1\right) ^{\left( \mu
-p-1\right) }\right) ^{\gamma }\right) ^{\frac{r}{\gamma }}\right] ^{\frac{1%
}{r}}\cdot \left( \sum_{s=a+\nu }^{b-\mu +p}\left| \Delta _{\ast }^{\mu
}f\left( s\right) \right| ^{\delta }\right) ^{\frac{1}{\delta }}.  \tag{33}
\label{r33}
\end{equation}
\end{theorem}

\begin{proof}
By (\ref{r31}) and discrete H\"{o}lder's inequality, we have 
\begin{equation*}
\left| \Delta ^{p}f\left( j\right) \right| \leq \frac{1}{\Gamma \left( \mu
-p\right) }\left( \sum_{s=a+\nu }^{j-\mu +p}\left( \left( j-s-1\right)
^{\left( \mu -p-1\right) }\right) ^{\gamma }\right) ^{\frac{1}{\gamma }}
\end{equation*}%
\begin{equation}
\cdot \left( \sum_{s=a+\nu }^{b-\mu +p}\left| \Delta _{\ast }^{\mu }f\left(
s\right) \right| ^{\delta }\right) ^{\frac{1}{\delta }},\text{ \ \ \ }%
\forall \text{ }j\in \lbrack a+m-p,b],  \tag{34}  \label{r34}
\end{equation}%
where $\gamma ,\delta >1:\frac{1}{\gamma }+\frac{1}{\delta }=1.$

Hence, by $r\geq 1$ we obtain 
\begin{equation*}
\left| \Delta ^{p}f\left( j\right) \right| ^{r}\leq \frac{1}{\left( \Gamma
\left( \mu -p\right) \right) ^{r}}\left( \sum_{s=a+\nu }^{j-\mu +p}\left(
\left( j-s-1\right) ^{\left( \mu -p-1\right) }\right) ^{\gamma }\right) ^{%
\frac{r}{\gamma }}
\end{equation*}%
\begin{equation}
\cdot \left( \sum_{s=a+\nu }^{b-\mu +p}\left| \Delta _{\ast }^{\mu }f\left(
s\right) \right| ^{\delta }\right) ^{\frac{r}{\delta }},\text{ \ \ \ }%
\forall \text{ }j\in \lbrack a+m-p,b].  \tag{35}  \label{r35}
\end{equation}%
Consequently we get 
\begin{equation*}
\sum_{j=a+m-p}^{b}\left| \Delta ^{p}f\left( j\right) \right| ^{r}\leq \frac{1%
}{\left( \Gamma \left( \mu -p\right) \right) ^{r}}\left[ \sum_{j=a+m-p}^{b}%
\left( \sum_{s=a+\nu }^{j-\mu +p}\left( \left( j-s-1\right) ^{\left( \mu
-p-1\right) }\right) ^{\gamma }\right) ^{\frac{r}{\gamma }}\right]
\end{equation*}%
\begin{equation}
\cdot \left( \sum_{s=a+\nu }^{b-\mu +p}\left| \Delta _{\ast }^{\mu }f\left(
s\right) \right| ^{\delta }\right) ^{\frac{r}{\delta }}.  \tag{36}
\label{r36}
\end{equation}%
The last proves the claim.
\end{proof}

We finish with the following discrete fractional average Sobolev type
inequali- ty.

\begin{theorem}
Let $0<\mu _{1}<\mu _{2}<...<\mu _{k}$; $m_{l}=\left\lceil \mu
_{l}\right\rceil $, $\nu _{l}=m_{l}-\mu _{l}$, $l=1,...,k,$ $k\in \mathbb{N}$%
. Assume that $\Delta ^{\tau }f\left( a\right) =0$, for $\tau
=0,1,...,m_{k}-1$: $f$ is defined on $\mathbb{N}_{a}$, $a\in \mathbb{Z}^{+}$%
. Let $r\geq 1$; $C_{l}\left( s\right) >0$ defined on $[a+\nu _{l},b-\mu
_{l}]$, $l=1,...,k$. Call

\noindent $B_{l}:=\sum_{s=a+\nu _{l}}^{b-\mu _{l}}C_{l}\left( s\right)
\left( \Delta _{\ast }^{\mu _{l}}f\left( s\right) \right) ^{2},$

\noindent $\delta ^{\ast }:=\underset{1\leq l\leq k}{\max }\left\{ \frac{1}{%
\left( \Gamma \left( \mu _{l}\right) \right) ^{2}}\left[
\sum_{j=a+m_{l}}^{b}\left( \sum_{s=a+\nu _{l}}^{j-\mu _{l}}\left( \left(
j-s-1\right) ^{\left( \mu _{l}-1\right) }\right) ^{2}\right) ^{\frac{r}{2}}%
\right] ^{\frac{2}{r}}\right\} $,

\noindent $\varrho ^{\ast }:=\underset{1\leq l\leq k}{\max }\left| \left( 
\frac{1}{C_{l}\left( s\right) }\right) \right| _{\infty ,[a+\nu _{l},b-\mu
_{l}]}.$

Then 
\begin{equation}
\left\Vert f\right\Vert _{r,[a+m_{k},b]}\leq \sqrt{\delta ^{\ast }\varrho
^{\ast }}\left( \frac{\sum_{l=1}^{k}B_{l}}{k}\right) ^{\frac{1}{2}}. 
\tag{37}  \label{r37}
\end{equation}
\end{theorem}

\begin{proof}
We see that also $\Delta ^{\tau }f\left( a\right) =0$, $\tau
=0,1,...,m_{l}-1 $, $l=1,...,k-1$. So the assumptions of Theorem \ref{t22}
are fulfilled for $f$ and fractional orders $\mu _{l}$, $l=1,...,k$. Thus by
choosing $p=0$ and $\gamma =\delta =2$ we apply (\ref{r33}), for $l=1,...,k$%
, to obtain 
\begin{equation*}
\left( \sum_{j=a+m_{l}}^{b}\left\vert f\left( j\right) \right\vert
^{r}\right) ^{\frac{1}{r}}\leq \frac{1}{\Gamma \left( \mu _{l}\right) }\left[
\sum_{j=a+m_{l}}^{b}\left( \sum_{s=a+\nu _{l}}^{j-\mu _{l}}\left( \left(
j-s-1\right) ^{\left( \mu _{l}-1\right) }\right) ^{2}\right) ^{\frac{r}{2}}%
\right] ^{\frac{1}{r}}
\end{equation*}%
\begin{equation}
\cdot \left( \sum_{s=a+\nu _{l}}^{b-\mu _{l}}\left( \Delta _{\ast }^{\mu
_{l}}f\left( s\right) \right) ^{2}\right) ^{\frac{1}{2}}.  \tag{38}
\label{r38}
\end{equation}%
Hence it holds 
\begin{equation*}
\left( \sum_{j=a+m_{l}}^{b}\left\vert f\left( j\right) \right\vert
^{r}\right) ^{\frac{2}{r}}\leq
\end{equation*}%
\begin{equation*}
\frac{1}{\left( \Gamma \left( \mu _{l}\right) \right) ^{2}}\left[
\sum_{j=a+m_{l}}^{b}\left( \sum_{s=a+\nu _{l}}^{j-\mu _{l}}\left( \left(
j-s-1\right) ^{\left( \mu _{l}-1\right) }\right) ^{2}\right) ^{\frac{r}{2}}%
\right] ^{\frac{2}{r}}\cdot \left( \sum_{s=a+\nu _{l}}^{b-\mu _{l}}\left(
\Delta _{\ast }^{\mu _{l}}f\left( s\right) \right) ^{2}\right)
\end{equation*}%
\begin{equation*}
\leq \delta ^{\ast }\left( \sum_{s=a+\nu _{l}}^{b-\mu _{l}}\left( \Delta
_{\ast }^{\mu _{l}}f\left( s\right) \right) ^{2}\right) =\delta ^{\ast
}\left( \sum_{s=a+\nu _{l}}^{b-\mu _{l}}\left( C_{l}\left( s\right) \right)
^{-1}\left( C_{l}\left( s\right) \right) \left( \Delta _{\ast }^{\mu
_{l}}f\left( s\right) \right) ^{2}\right)
\end{equation*}%
\begin{equation}
\leq \delta ^{\ast }\rho ^{\ast }\left( \sum_{s=a+\nu _{l}}^{b-\mu
_{l}}C_{l}\left( s\right) \left( \Delta _{\ast }^{\mu _{l}}f\left( s\right)
\right) ^{2}\right) .  \tag{39}  \label{r39}
\end{equation}

That is 
\begin{equation*}
\left( \sum_{j=a+m_{k}}^{b}\left\vert f\left( j\right) \right\vert
^{r}\right) ^{\frac{2}{r}}\leq \left( \sum_{j=a+m_{l}}^{b}\left\vert f\left(
j\right) \right\vert ^{r}\right) ^{\frac{2}{r}}\leq
\end{equation*}%
\begin{equation}
\delta ^{\ast }\rho ^{\ast }\left( \sum_{s=a+\nu _{l}}^{b-\mu
_{l}}C_{l}\left( s\right) \left( \Delta _{\ast }^{\mu _{l}}f\left( s\right)
\right) ^{2}\right) =\delta ^{\ast }\rho ^{\ast }B_{l},\text{ \ \ \ for }%
l=1,...,k.  \tag{40}  \label{r40}
\end{equation}

Hence 
\begin{equation}
\left\| f\right\| _{r,[a+m_{k},b]}^{2}\leq \delta ^{\ast }\rho ^{\ast
}\left( \frac{\sum_{l=1}^{k}B_{l}}{k}\right) ,  \tag{41}  \label{r41}
\end{equation}%
proving the claim.
\end{proof}

\end{document}